\documentclass{article}
\usepackage{amsmath,amsfonts,latexsym,amssymb,euscript,amsthm}

\newcommand{\R}{\mathbb{R}}
\newcommand{\Om}{\Omega}

\renewcommand{\d}{\partial}
\newcommand{\vphi}{\varphi}
\newcommand{\ep}{\varepsilon}
\newcommand{\intl}{\int\limits}
\renewcommand{\div}{\mathop{\mathrm{div }}}
\newcommand{\rot}{\mathop{\mathrm{rot }}}

\newtheorem{theorem}{Theorem}[section]
\newtheorem{definition}{Definition}[section]
\newtheorem{lemma}{Lemma}[section]

\renewcommand{\l}{\left}
\renewcommand{\r}{\right}

\newcommand{\Q}{Q^+}
\newcommand{\B}{B^+}
\newcommand{\al}{\alpha}
\newcommand{\WO}{\stackrel{0}{W^1_2}}
\newcommand{\gr}{\nabla}

\newcommand{\lap}{\Delta}
\renewcommand{\div}{\mathop{\mathrm{div }}}
\newcommand{\Ga}{\Gamma}
\newcommand{\Bdva}{B^{(2)}}
\newcommand{\Qdva}{Q^{(2)}}

\newcommand{\ph}{\varphi}

\renewcommand{\L}{\EuScript{L}}

\newcommand{\A}{\EuScript{A}}

\newcommand{\cd}{\partial}

\newcommand{\tH}{\tilde{H}}
\newcommand{\tv}{\tilde{v}}
\newcommand{\tpsi}{\tilde{\psi}}

\newcommand{\be}{\beta}
\newcommand{\ga}{\gamma}

\newcommand{\hi}{h^{(1)}}
\newcommand{\hii}{h^{(2)}}
\newcommand{\hiii}{h^{(3)}}

\begin{document}

\title{On the Regularity of Weak Solutions to the Magneto Hydrodynamics System near the curved part of the boundary}
\author{V.~Vialov
\footnote{This research is supported by the Chebyshev Laboratory
(Department of Mathematics and Mechanics, St.-Petersburg State University)
under RF government grant 11.G34.31.0026}}

\date{\today}

\maketitle
\abstract{ We prove a sufficient conditions of local regularity of suitable weak solutions to the MHD system for the point from $C^3$-smooth part of the boundary. Our
conditions are the generalizing of the Caffarelli--Kohn--Nirenberg theorem for Navier-Stokes equations.}

\section{Introduction}

Assume $\Om\subset \mathbb R^3$ is a $C^3$-- smooth bounded domain
and $Q_T=\Om\times (0,T)$. In this paper we investigate the boundary
regularity of solutions to the principal system of
magnetohydrodynamics (the MHD equations):
\begin{equation}
\left.
\begin{array}c
\partial_t v  + (v\cdot \nabla )v  - \Delta v +  \nabla
p = \rot H \times H
\\
 \div v =0
\end{array} \right\}
\quad\mbox{in } Q_T, \label{MHD_NSE}
\end{equation}
\begin{equation}
\left.
\begin{array}c
\partial_t H    +\rot \rot H = \rot (v\times H)
\\ \div H =0
\end{array} \right\} \quad \mbox{in  }Q_T.
\label{MHD_Magnetic}
\end{equation}
Here  unknowns
are  the velocity field $v:Q_T\to \mathbb R^3$, pressure
$p:Q_T\to \mathbb R$, and the magnetic field $H:Q_T\to \mathbb
R^3$. We impose on $v$ and $H$ the boundary conditions:
\begin{equation}
v|_{\cd\Om\times(0,T)}=0, \quad H_{\nu}|_{\cd\Om\times(0,T)}=0, \quad (\rot H)_{\tau}|_{\cd\Om\times(0,T)}=0,
\label{BC}
\end{equation}
Here by $\nu$ we denote the outer normal to $\cd\Om$ and $H_{\nu}= H\cdot \nu$, $(\rot H)_\tau=\rot H- \nu(\rot H\cdot \nu)$.
These conditions correspond to the case of liquid flowing in the area bounded by ideal conductor.


\begin{definition}
\label{sws_def}
Assume $\Ga\subset \cd\Om$.
 The functions $(v,H,p)$ are called a {\it boundary suitable
 weak solution } to the system  (\ref{MHD_NSE}), (\ref{MHD_Magnetic})
near $\Ga_T\equiv \Ga\times (0,T)$ if
\begin{itemize}
\item[1)] \ $v\in L_{2,\infty}(Q_T)\cap W^{1,0}_2(Q_T)\cap W^{2,1}_{\frac 98, \frac 32}(Q_T)$,

$H \in L_{2,\infty}(Q_T)\cap W^{1,0}_2(Q_T)$,
\item[2)] \ $p\in L_{\frac 32}(Q_T)\cap W^{1,0}_{\frac 98, \frac 32}(Q_T)$,
\item[3)] \ $\div v=0$, \ $\div H=0$ \ a.e. in \ $Q_T$,
\item[4)] \ $v|_{\cd \Om}=0$, \ $H_\nu|_{\cd\Om}=0$ \ in the sense of traces,
\item[5)] \ for any $w\in L_2(\Om)$ the functions
$$
t\mapsto \int\limits_{\Om} v(x,t)\cdot w(x)~dx \qquad\mbox{and}\qquad t\mapsto \int\limits_{\Om} H(x,t)\cdot w(x)~dx
$$
are continuous,

\item[6)] \ $(v,H)$ satisfy the following integral identities: for any $t\in [0,T]$
$$
\begin{array}c
\int\limits_{\Om} ~v(x,t)\cdot \eta(x,t)~dx -\int\limits_{\Om} ~v_0(x)\cdot \eta(x,0)~dx \ + \\ + \
\int\limits_0^t \int\limits_\Om ~ \Big( - v\cdot \cd_t \eta + (\nabla v  - v\otimes v + H \otimes H) : \nabla \eta - (p+\frac 12 |H|^2)\div \eta \Big)~dxdt \ = \ 0,
\label{isl1}
\end{array}
$$
for all \ $\eta \in W^{1,1}_{\frac 52}(Q_t)$ \ such that   \ $\eta|_{\cd\Om\times (0,t)}=0$,
$$
\begin{array}c
 \int\limits_{\Om}~ H(x,t)\cdot \psi (x,t)~dx  \ -  \  \int\limits_{\Om}~ H_0(x)\cdot \psi (x,0)~dx \ + \\ + \
\int\limits_0^t \int\limits_\Om ~ \Big( - H\cdot \cd_t \psi  + \rot H\cdot \rot \psi - (v\times H) \cdot \rot \psi \Big)~dxdt  \ = 0,
\label{isl2}
\end{array}
$$
for all \ $\psi \in W^{1,1}_{\frac 52}(Q_t)$ \ such that \  $\psi_\nu |_{\cd\Om\times (0,t)}=0$.

\item[7)] For every $z_0=(x_0,t_0)\in \Ga_T$
such that  $ \Om_R(x_0)\times (t_0-R^2, t_0)\subset Q_T$ where $\Om(x_0,R) \equiv \Om \cap B(x_0,R)$ and for any $\zeta\in C_0^\infty(B_R(x_0)\times(t_0-R^2,t_0])$
such that $\left.\frac{\cd\zeta}{\cd \nu}\right|_{\cd\Om}=0$ the following ``local energy inequality near $\Ga_T$'' holds:
\begin{equation}
\begin{array}c
\sup\limits_{t\in (t_0-R^2,t_0)}\int\limits_{\Om_R(x_0)}  \zeta \Big( |v|^2+|H|^2\Big)~dx \ + \\
+ \ 2 \int\limits_{t_0-R^2}^{t_0}\int\limits_{\Om_R(x_0)} \zeta \Big( |\nabla v|^2+|\rot H|^2\Big)~dxdt \ \le \\
\le\  \int\limits_{t_0-R^2}^{t_0}\int\limits_{\Om_R(x_0)}  \Big(| v|^2 +|H|^2\Big) (\cd_t \zeta + \Delta \zeta )~dxdt \ + \\ + \
\int\limits_{t_0-R^2}^{t_0}\int\limits_{\Om_R(x_0)} \Big ( | v|^2+ 2\bar p \Big ) v\cdot\nabla \zeta~dxdt \ + \\
- \ 2 \int\limits_{t_0-R^2}^{t_0}\int\limits_{\Om_R(x_0)}   (H\otimes H):  \nabla^2 \zeta  ~dxdt \ + \\ + \
2 \int\limits_{t_0-R^2}^{t_0}\int\limits_{\Om_R(x_0)}  (v\times  H)(\nabla \zeta\times H)~dxdt
\end{array}
\label{LEI}
\end{equation}
\end{itemize}
\end{definition}

\noindent
Here  $L_{s,l}(Q_T)$ is the anisotropic Lebesgue space equipped with the norm
$$
\|f\|_{L_{s,l}(Q_T)}:=
\Big(\int_0^T\Big(\int_\Om |f(x,t)|^s~dx\Big)^{l/s}dt\Big)^{1/l} ,
$$
and we use the following notation for the functional spaces:
$$
\gathered
W^{1,0}_{s,l}(Q_T)\equiv L_l(0,T; W^1_s(\Om))= \{ \ u\in L_{s,l}(Q_T): ~\nabla u \in L_{s,l}(Q_T) \ \},
\\
W^{2,1}_{s,l}(Q_T) = \{ \ u\in W^{1,0}_{s,l}(Q_T): ~\nabla^2 u, \ \cd_t u \in L_{s,l}(Q_T) \ \},
\\
\overset{\circ}{W}{^1_s}(\Om)=\{ \ u\in W^1_s(\Om):~ u|_{\cd\Om}=0 \ \},
\endgathered
$$
and the following notation for the norms:
$$
\gathered
\| u \|_{W^{1,0}_{s,l}(Q_T)}= \| u \|_{L_{s,l}(Q_T)}+ \|\nabla u\|_{L_{s,l}(Q_T)},  \\
\| u \|_{W^{2,1}_{s,l}(Q_T)}= \| u \|_{W^{1,0}_{s,l}(Q_T)}+ \| \nabla^2 u \|_{L_{s,l}(Q_T)}+\|\cd_t u\|_{L_{s,l}(Q_T)},  \\
\endgathered
$$

\begin{theorem} \label{Existense_Suitable}
For any sufficiently smooth divergent-free \ $v_0$,  $H_0$ satisfying (\ref{BC})  there exists at least one boundary suitable weak solution near $\cd\Om\times (0,T)$ which satisfies the initial conditions :
$$\| v(\cdot, t) -v_0(\cdot)\|_{L_2(\Om)}\to 0, \quad \| H(\cdot, t) -H_0(\cdot)\|_{L_2(\Om)}\to 0 \quad\mbox{as} \quad t\to+0,$$
and additionally satisfies the global energy inequality
$$
\gathered
\| v\|_{L_{2,\infty}(Q_T)}+ \| H\|_{L_{2,\infty}(Q_T)} + \| \nabla v\|_{L_{2}(Q_T)} + \| \rot H\|_{L_{2}(Q_T)} \le \\ \le \| v_0\|_{L_2(\Om)}+ \| H_0\|_{L_2(\Om)}
\endgathered
$$
\end{theorem}

The global existence of weak solutions to the MHD equations (\ref{MHD_NSE}) --- (\ref{BC}) was established originally  in \cite{LadSol}.
The proof of Theorem \ref{Existense_Suitable} can be found in \cite{VyaShi}.

\noindent
{\bf Notations} In this paper we will use the following notations
$$B(x_0,R) = \l\{ x \in \R^3 : |x-x_0| < R \r\},$$
$$\Bdva(x'_0,R) = \l\{ x' \in \R^3 : |x'-x'_0| < R \r\},$$
$$\Bdva(R) = \Bdva(0,R), \quad \Bdva = \Bdva(1)$$
$$\B_0(R) = \l\{ x \in B(0,R) : x_3 > 0 \r\}, \quad \Q_0(R) = \B_0(R) \times (-R^2,0)$$


\section{Main Results}
\setcounter{equation}{0}

\noindent
{\bf Main condition on $\d \Om$}. {\it There exist positive numbers $\mu$ and $R_0$ depending only on $\d \Om$ such that
for each point $x_0 \in \d \Om$ we can choose a Cartesian coordinate system $\{y_i\}^3_{i=1}$ associated to the origin $x_0$, and some function $\ph_{x_0} \in C^3(\Bdva(R_0))$ such that
$$\Om(x_0,R_0) \equiv \Om \cap B(x_0,R_0) = \{ y \in B(R_0) : y_3 > \ph_{x_0}(y_1,y_2) \},$$
and
\begin{equation}
\vphi_{x_0}(0) = 0, \quad \gr \vphi_{x_0} (0) = 0, \quad \| \vphi_{x_0} \|_{W^2_{\infty}} \leq \mu.
\label{phi_cond}
\end{equation}
}


The main results of the present paper are the following theorems on boundary regularity of suitable weak solutions of MHD system
\begin{theorem}
\label{CKN_theorem_Omega}
Assume that $\d\Om$ satisfies Main Condition. Then for any  $K>0$ there exists $\ep_0(K)>0$ with the following property.
Assume $(v,H,p)$ is a boundary suitable weak solution in $Q_T$ and $z_0=(x_0,t_0)\in \cd\Om\times (0,T)$.
If
\begin{equation}
\limsup\limits_{r\to 0}\Big(~\frac 1r \int\limits_{t_0-r^2}^{t_0} |\nabla H|^2~dxdt ~\Big)^{1/2}\ < \ K
\end{equation}
and
\begin{equation}
\limsup\limits_{r\to 0} \Big(~\frac 1r \int\limits_{t_0-r^2}^{t_0} |\nabla v|^2~dxdt~\Big)^{1/2} \ < \ \ep_0,
\end{equation}
then the functions $v$ and $H$ are H\" older continuous in some neighborhood of $z_0$.
\end{theorem}

\begin{theorem}\label{Partial_Regularity}
Assume that $\d\Om$ satisfies Main Condition and $(v,H,p)$ is a boundary suitable weak solution in $Q_T$. Then there exists a closed set $\Sigma\subset \d\Om \times (0,T]$ such that
for any $z_0\in (\Ga\setminus \Sigma)\times (0,T]$ the functions $(v,H)$ are H\" older continuous in some neighborhood of $z_0$,
\begin{equation}
\mathcal P^1(\Sigma) \ = \ 0,
\label{Sigma_estimate}
\end{equation}
where $\mathcal P^1(\Sigma)$ is the one-dimensional parabolic Hausdorff measure of $\Sigma $.
\end{theorem}

Our Theorem \ref{CKN_theorem_Omega} presents for the MHD equations a result which is a boundary analogue of  the famous Caffarelli--Kohn--Nirenberg (CKN)
theorem  for the Navier-Stokes system, see \cite{CKN}, see also \cite{Lin}.
The boundary regularity of solutions to the Navier-Stokes equations was originally investigated  by G.~Seregin in \cite{Seregin_JMFM} and \cite{Seregin_Aa}  in the case of a plane part of the boundary
and by G.~Seregin, T.~Shilkin, and V.~Solonnikov in \cite{SSS} in the case of a curved boundary.

The internal partial regularity of  solutions to the MHD system was originally proved by C.~He and Z.~Xin in \cite{China}, see also \cite{Vya}, \cite{Vya1}. The local regularity near the plan
part of the boundary for MHD equations was investigated in \cite{VyaShi} (see also \cite{Vya2011}). In the case of boundary regularity due to boundary conditions on the derivatives of $H$
there will be some problems if try to directly generalize the approach used for Navier-Stokes equations. To solve this problem we will consider the equation \eqref{MHD_Magnetic} as a linear heat equation on $H$
and this gives us some additional estimates. This idea was originally proposed in \cite{China} to obtain regularity theorems with just boundedness conditions instead of smallness on magnetic component  the internal case.
In the present paper corresponding estimates are proved in the sections \ref{Heat_equation_estimate_section} and \ref{H_estimate_section}.

Note that using the methods of our paper one can prove various $\ep$--regularity conditions involving various scale--invariant  functionals
(such it was done for the plane part of the boundary in \cite{Vya2011}, see also \cite{Mih}).
In the present paper we concentrate on the condition of the theorem \ref{CKN_theorem_Omega} as this condition provides the optimal estimate of the Hausdorff measure of the singular set $\Sigma$
in Theorem \ref{Partial_Regularity}. Note that we also have \eqref{Sigma_estimate} in the internal case, so combining these two estimates we will obtain \eqref{Sigma_estimate} for the set of singular points
in any bounded domain $\Omega$ with $C^3$ smooth boundary.

Our paper is organized as follows:
in Section \ref{Symmetry_section} using symmetries of \eqref{MHD_NSE}-\eqref{MHD_Magnetic} we present more convenient statement of Theorem \ref{CKN_theorem_Omega}.
In Section \ref{Perturbed_Stokes_section} we describe coercive estimates for solutions of the Stokes equations near the boundary.
Section \ref{Main_criterion_section} contains the proof of the Decay Lemma and the sketch of the proof of Theorem \ref{Fixed_r}.
Sections \ref{Heat_equation_estimate_section} and \ref{H_estimate_section} is concerned with the estimate of some Morrey functional for weak solutions to the heat equation near the boundary.
These estimates together with the estimates of the scale invariant energy functionals obtained in Section \ref{EnergyEstimates_Section} turn to be crucial for the prove of theorem \ref{CKN_theorem_Omega}
presented in section \ref{Final_section}.


\section{Symmetry group of MHD system and new statement of the main results}
\label{Symmetry_section}
\setcounter{equation}{0}

The solutions of MHD system (\ref{MHD_NSE}), (\ref{MHD_Magnetic}) have the same set of symmetries as the Navier-Stokes equations i.e.
they are invariant under translations, rotations and scaling
\begin{equation}
\begin{aligned}
& v_R(y,s) = R v(R y, R^2 s),\\
& H_{R}(y,s) = R H(R y, R^2 s),\\
& p_{R}(y,s) = R^2 p(R y, R^2 s).
\label{scaling}
\end{aligned}
\end{equation}
So we can consider that in the statement of theorem \ref{CKN_theorem_Omega}  $z_0 = 0$ and the boundary of the domain is described by function $\ph$ satisfying \eqref{phi_cond}.
Also it will be convenient to consider the function $\ph$ as a part of the problem data
and deal with the "local version" of suitable weak solution in parabolic cylinder
$$\Q(R) \equiv (-R^2,0) \times \B(R),$$
where
$$\B(R) \equiv \{ x \in B(R) : x_3 > \ph (x_1,x_2) \}.$$

\begin{definition}
Let $R>0$ and $\vphi \in C^2(\Bdva(R))$ satisfies \eqref{phi_cond}.
The functions $(v,H,p,\vphi)$ are called a boundary suitable
 weak solution to the system  (\ref{MHD_NSE}), (\ref{MHD_Magnetic})
in $\Q(R)$ if there is a domain $\Om$ such that $\Gamma = \l\{ x_3 = \ph(x_1,x_2) \r\} \subseteq \d \Om$ and functions $v$, $p$, $H$ can be extended up to
suitable weak solution near $\Gamma$.
\end{definition}


Then theorem \ref{CKN_theorem_Omega} can be formulated by following way

\begin{theorem}
\label{CKN_theorem}
For any  $K>0$ there exists $\ep_0(K)>0$ with the following property.
Assume $(v,H,p,R)$ is a boundary suitable weak solution in $\Q(R)$ for some $R>0$.
If
\begin{equation}
\limsup\limits_{r\to 0}\Big(~\frac 1r \int\limits_{\Q(r)} |\nabla H|^2~dxdt ~\Big)^{1/2}\ < \ K
\label{ep-regularity-1}
\end{equation}
and
\begin{equation}
\limsup\limits_{r\to 0} \Big(~\frac 1r \int\limits_{\Q(r)} |\nabla v|^2~dxdt~\Big)^{1/2} \ < \ \ep_0,
\label{ep-regularity-2}
\end{equation}
then there exists $\rho_*>0$ such that the functions $v$ and $H$ are H\" older continuous on the closure of $ Q^+(\rho_*)$.
\end{theorem}

To prove this theorem we will generalize the approach introduced in \cite{VyaShi} in the case of the plane part of the boundary. The first step is the following theorem

\begin{theorem}\label{Fixed_r}
There exists an absolute constant $\ep_*>0$ with the following property. Assume $(v,H,p,\vphi)$ is a boundary suitable weak solution in $\Q(R)$
If there exists $0 < r_0 <R$ such
$$
\frac 1{r_0^2} \int\limits_{Q^+(r_0)} \Big(~ |v|^3+|H|^3+|p|^{\frac 32}~\Big) dxdt \  < \ \ep_*
$$
and
\begin{equation}
\| \ph \|_{C^2(\Bdva(r_0))} < \ep_*,
\label{ssl1}
\end{equation}
then the functions $v$ and $H$ are H\" older continuous on the closure of $ Q^+(\frac {r_0}2)$.
\end{theorem}

Note, that \eqref{ssl1} is just the condition on the smallness of $r_0$. Indeed
if $(v,p,H,\vphi)$ are the suitable weak solution in $\Q(R)$, then if we apply the scaling transformations \eqref{scaling},
then $(v_R,p_R,H_R,\ph_R)$ where
$$
\vphi_R  = \frac1{R}\vphi
$$
will be the solution in $\Q(1)$ and from Taylor formula we have
$$ \| \frac{\vphi}{R} \|_{C^2 (\Bdva)} \leq R \| \vphi \|_{C^2}.$$


\section{Estimates for perturbed Stokes system}
\label{Perturbed_Stokes_section}
\setcounter{equation}{0}

In this section we describe coercive estimates for linearisation of the \eqref{MHD_NSE}. We start from the Stokes problem in $\Q$ with
some $\ph \in C^3(\Bdva)$
\begin{equation}
\begin{aligned}
& \d_t u - \Delta u + \nabla p = f\\
& \nabla \cdot u = 0 \\
& u|_{x_3=\ph(x_1,x_2)}  = 0,
\end{aligned}
\qquad \text{ in } \Q.
\label{pssl1}
\end{equation}
and define new coordinates $\{ y_i \}_{i=1}^3$ connected with the original ones via formula
\begin{equation}
x = e (y) \equiv \begin{pmatrix} y_1\\ y_2 \\ y_3 + \vphi(y_1,y_2) \end{pmatrix}.
\label{change_variables}
\end{equation}
Denote by $\L$ the Jacobi matrix of the map $x = e(y)$ i.e.
\begin{equation}
\L =
\begin{pmatrix}
1 & 0 & 0\\
0 & 1 & 0\\
\vphi_{,1} & \vphi_{,2} & 1
\end{pmatrix},
\label{Jmatrix}
\end{equation}
and by $\tilde \nabla_\ph$ and $\tilde \Delta_\ph $ the following differential operators with variable coefficients:
\begin{equation}
\gathered
(\tilde \nabla_\ph p)_{i} = p_{,k}{e_{k,i}}, \\ (\tilde \Delta_\ph v)_i = v_{i,jl}e_{j,k}e_{l,k} +v_{i,j}e_{j,kk}, \\
\tilde\nabla_\ph \cdot v  = v_{i,k}{e_{k,i}}
\endgathered
\label{Hat operators-1}
\end{equation}
Also we have the relation
\begin{equation}
\rot\nolimits_x H = \L \rot\nolimits_y (\L^T \tH),
\label{rot_tr}
\end{equation}
here we have used that $\det \L = 1$.

Then the system \eqref{pssl1} transforms into so called perturbed Stokes system
\begin{equation}
\begin{aligned}
& \d_t u - \tilde\Delta_{\ph} u + \tilde\nabla_{\ph} p = f\\
& \tilde\nabla_{\ph} \cdot u = 0 \\
& u|_{y_3=0}  = 0.
\end{aligned}
\label{zerodiv_Stokes}
\end{equation}
Without loss of generality we can consider \eqref{zerodiv_Stokes} in $\Q_0$.

We recall that the function $\vphi$ satisfy the following relations
\begin{equation}
\vphi (0) = 0, \qquad \nabla \vphi (0) = 0, \qquad \| \vphi \|_{W^2_{\infty}} < \mu,
\label{phi_conditions}
\end{equation}
and we can think, that the constant $\mu$ is sufficiently small.

For our further arguments we will need two lemmas about the solutions of this problem.
We start from consideration of initial-boundary problem for perturbed Stokes system \eqref{zerodiv_Stokes} with homogenous data:
\begin{equation}
u|_{\d \B}  = 0, \quad u|_{t=-1} = 0
\label{nonzerodiv_Stokes_BC}
\end{equation}

\begin{lemma}
Assume that $s,l \in (1,\infty)$. There is $\mu_* \in (0,\frac1{100})$ if $\vphi$ satisfies \eqref{phi_conditions} for some $\mu < \mu_*$,
then there is a pair of functions $(u,p)$ such that
$$ u \in W^{2,1}_{s,l}(\Q), \qquad \nabla p \in L_{s,l}(\Q),$$
$(v,p)$ satisfy \eqref{zerodiv_Stokes} with initial-boundary conditions \eqref{nonzerodiv_Stokes_BC} and the following estimate holds
\begin{equation}
\| v \|_{W^{2,1}_{s,l}(\Q)} + \| \nabla p \|_{L_{s,l}(\Q)} \leq C_*  \| f \|_{L_{s,l}(\Q)}.
\label{linest1}
\end{equation}
here $C_*$ is an absolute constant depending only on the dimension.
\label{nonzerodivStokes_Lemma}
\end{lemma}

The second lemma is the coercive estimate for the solution of \eqref{zerodiv_Stokes}.

\begin{lemma}
Assume that $s,m,l \in (1,\infty)$, $m \geq s$. Assume that $\vphi \in C^3(\Bdva)$ satisfies \eqref{phi_conditions} for some $\mu < \mu_*$.
Then for every functions $u \in W^{2,1}_{s,l}(\Q)$, $\nabla p \in L_{s,l}(\Q)$ and $f \in L_{m,l}(\Q)$ satisfying \eqref{zerodiv_Stokes} we have
$u \in W^{2,1}_{m,l}(\Q(\frac12))$, $\nabla p \in L_{m,l}(\Q(\frac(1/2))$ and the following estimate holds
\begin{equation}
\gathered
\| v \|_{W^{2,1}_{m,l}(\Q(\frac12))} + \| \nabla p \|_{L_{m,l}(\Q(\frac12))} \leq\\
\leq C \l( \| f \|_{L_{s,l}(\Q)} + \|\nabla v \|_{L_{s,l}(\Q)} + \| p - p_0 \|_{L_{s,l}(\Q)} \r)
\endgathered
\label{linest2}
\end{equation}
for some absolute constant $C$ and arbitrary function $p_0 = p_0 (t)$, $p \in L_l(-1,0)$.
\label{zerodivStokes_Lemma}
\end{lemma}


\section{The main criterion of $\ep$-regularity}
\label{Main_criterion_section}
\setcounter{equation}{0}

In this section we will prove the theorem \ref{Fixed_r}. Let $(v,p,H,\ph)$ are the suitable weak solution in $\Q$, then we will use the following notations
$$
\gathered
Y_\tau (v) \ =  \l(\frac{1}{|Q(\tau)|} \int\limits_{Q^+(\tau)} ~|v|^3~dxdt\ \r)^{1/3}, \\
\tilde Y_\tau (H) \ = \l(\frac{1}{|Q(\tau)|} \int\limits_{Q^+(\tau)} ~|H - b_R(H)|^3~dxdt\ \r)^{1/3}, \\
\hat Y_\tau (p)\ = \
\tau ~\l(\frac{1}{|Q(\tau)|}\int\limits_{Q^+(\tau)} ~|p-[p]_{B^+(\tau)}|^{3/2}~dxdt\ \r)^{2/3},\\
Y_\tau (v,p,H) = Y_\tau(v) + \tilde Y_\tau (H) + \hat Y_\tau (p),\\
b_\tau(H) = \L \begin{pmatrix} (h_1)_\tau \\ (h_2)_\tau \\ 0 \end{pmatrix},
\endgathered
$$
here $h = \L (H \circ e)$, where $e$ and $\L$ defined by \eqref{change_variables} and \eqref{Jmatrix},
$$
(h_i)_\tau = \frac1{|\Q(\tau)|} \intl_{e(\Q(\tau))} h_i(z)\,dz
$$
and
$$
[p]_{\B(\tau)}(t) = \frac1{|\B(\tau)|} \intl_{\B(\tau)} p(x,t)\,dt.
$$
Note, that $\tilde Y_\tau(H)$ is equivalent to norm $H$ in Morrey space
with taking into account boundary conditions. It can be easily seen if we replace $H$ by $\L h$ and use, that matrix $\L$ is close to
the identity  one.


We start our considerations from the modification of local energy inequality.

\begin{lemma}
Assume  $(v,H,p,\ph)$ is a boundary suitable weak solution satisfied the MHD equations in $\Q(R)$.
Let $\zeta\in C_0^\infty(B\times (-R^2,0])$ be a cut-off function such that $\frac{\d\zeta}{\d \nu}|_{x_3 = \ph(x_1,x_2)}=0$.
Assume $b\in \mathbb R^3$ is an arbitrary constant vector of the form $b=(b_1, b_2, 0)$. Then the following inequality holds
\begin{equation}
\gathered
\intl_{\B(R)} \zeta \l( |v|^2 + | \bar{H}|^2 \r)\,dx \\
+ 2 \intl_{\B(R)\times (t_0,t)} \zeta\l( |\nabla v| + |\rot \bar{H}|^2\r)\,dz\\
\leq \intl_{\B(R)\times (t_0,t)} (\d_t \zeta + \lap \zeta) (|v|^2 + |\bar{H}|^2)\,dz\\
+\intl_{\B(R)\times (t_0,t)} \Big ( | v|^2+ 2\bar p \Big ) v\cdot\nabla \zeta\,dz \\
- \ 2 \intl_{\B(R)\times (t_0,t)}   (H\otimes \bar{H}):  \nabla^2 \zeta \,dz \ + \\ + \
2 \intl_{\B(R)\times (t_0,t)}  (v\times  H)(\nabla \zeta\times \bar{H})\,dz\\
- \intl_{\B(R)\times (t_0,t)} \l[ |\rot \L b|^2 \zeta + \nabla \zeta \cdot \nabla (|\L b|^2) - (v\times H) \cdot \rot(\L b) \zeta\r]\,dz,
\endgathered
\label{mlei1}
\end{equation}
where $\bar{H} = H - \L b$.
\end{lemma}

{\bf Proof.}
We use \eqref{LEI} and transform the remaining terms. Via integration by parts formula
\begin{equation}
\gathered
\intl_{\B(R)\times (t_0,t)} (\d_t \zeta + \lap \zeta) |\L b|^2\,dz \\
= \intl_{\B(R)} \zeta |\L b|^2 \,dx -
\intl_{\B(R)\times (t_0,t)} \nabla \zeta \cdot \nabla (|\L b|^2)\,dz
\endgathered
\label{mlei2}
\end{equation}
Also we have
\begin{equation}
\gathered
\intl_{\B(R)\times (t_0,t)} \d_t \zeta H \cdot \L b \,dz = \intl_{\B(R)} \zeta H \cdot \L b \,dx +\\
+ \intl_{\B(R)\times (t_0,t)} (-\d_t H \cdot \L b)\,dz.
\endgathered
\label{mlei3}
\end{equation}
Now we consider two terms: one from left hand side another from right hand side of \eqref{mlei1}
\begin{equation}
\gathered
\intl_{\B(R)\times (t_0,t)} \lap \zeta H \cdot \L b\,dz -
\intl_{\B(R)\times (t_0,t)} \rot H \cdot \rot (\L b)\,dz =\\
= - \intl_{\B(R)\times (t_0,t)} \l[ \rot H \cdot \l(\nabla \zeta \times (\L b) \r) + \rot H \cdot \rot (\L b) \zeta \r.\\
\l.- H \otimes (\L b) : \nabla^2 \zeta\r]\,dz\\
-\intl_{\B(R)\times (t_0,t)} H \otimes (\L b) : \nabla^2 \zeta\,dz=\\
= -\intl_{\B(R)\times (t_0,t)} H \otimes (\L b) : \nabla^2 \zeta\,dz
- \intl_{\B(R)\times (t_0,t)} \rot \rot H \cdot (\L b \zeta)\,dz.
\endgathered
\label{mlei4}
\end{equation}

The last terms of \eqref{mlei3} and \eqref{mlei4} can be modified using equation \eqref{MHD_Magnetic}
\begin{equation}
\gathered
\intl_{\B(R)\times (t_0,t)} (\d_t H \cdot \L b + \rot \rot H \cdot (\L b \zeta))\,dz=\\
= \intl_{\B(R)\times (t_0,t)} \rot (v \times H) \cdot (\L b\zeta)\,dz=\\
= \intl_{\B(R)\times (t_0,t)} (v \times H) \cdot \rot (\L b) \zeta -\\
- \intl_{\B(R)\times (t_0,t)} (v \times H) \cdot (\nabla \zeta \times (\L b))\,dz.
\endgathered
\label{mlei5}
\end{equation}
Combining \eqref{mlei2}-\eqref{mlei5} we obtain the statement of lemma.

\qed


\begin{lemma} \label{Decay estimate theorem}
There exists an absolute constant $\ep_*>0$ such that for any $M>0$ there exists $C_*=C_*(M)$ with the following properties.
For any boundary suitable weak solution $(v,H,p,\vphi)$ of the MHD system (\ref{MHD_NSE}), (\ref{MHD_Magnetic}) near the boundary in $\Q(1)$
the following implication holds:
if
$$
Y_1(v,H,p) + \| \vphi \|_{C^2} \ < \ \ep_0,
$$
and
$$
b_R(H) \le \  M
$$
then
\begin{equation}
Y_\tau (v,H,p) \ \le \ C_* ~\tau^{1/3}~Y_1 (v,H,p)
\label{Decay_Estimate}
\end{equation}
\end{lemma}

{\bf Proof.}
Arguing by contradiction we assume there exists a sequence of numbers $\ep_m\to 0$,
and a sequence of boundary suitable weak solutions $(v^m,H^m, p^m, \vphi^m)$ such that
$$
Y_1(v^m,H^m,p^m) + \| \vphi^m \|_{C^2} \ < \ \ep_m \ \to \ 0,
$$
and
$$
Y_\tau (v^m, H^m, p^m)\ \ge \  C_* \tau^{1/3} \delta_m
$$
here $\delta_m = Y_1(v^m,H^m,p^m)$.

Let us introduce functions
$$
\gathered
u^m(y, s) \ = \  \frac{1}{\delta_m}~ v^m \circ e_m, \\
q^m(y, s) \ = \ \frac{1}{\delta_m}~ \Big(p^m(x, t) - [p^m]_{B^+}(t)\Big) \circ e_m, \\
h^m(y, s) \ = \  \frac{1}{\delta_m}~ \Big(H^m(x, t) - b_1(H^m)\Big) \circ e^m,\\
\endgathered
$$
here $e^m$ denotes the map \eqref{change_variables} corresponding to $\vphi^m$.
Then
\begin{equation}
Y_1(u^m,h^m,q^m)\ = \  1,
\qquad
Y_\tau (u^m, h^m, q^m )\ge C_* \tau^{1/3}
\label{Contradiction_asumptions}
\end{equation}
and $(u^m, h^m, q^m)$ satisfy the following equations in $\mathcal D'(Q^+)$
\begin{equation}
\gathered
\partial_t u^m + \delta_m (u^m\cdot \tilde\nabla_{\vphi^m} )u^m  - \tilde\Delta_{\vphi^m} u^m +  \tilde\nabla_{\vphi^m} q^m =\\
= \L_m \rot \L_m^{-T} h^m \times ( \delta_m h^m +  b_1(H^m))
\\
 \L_m^{-T} \cdot \nabla u^m =0
\endgathered  \label{MHD_NS_m}
\end{equation}
\begin{equation}
\gathered
\partial_t h^m     -  \tilde\Delta_{\vphi^m} h^m = \L_m \rot \L_m^{T} \big(u^m\times (\delta_m h^m+ (H^m))\big)
\\ \L_m^{-T} \cdot \nabla h^m =0
\endgathered
\label{MHD_Magnetic_m}
\end{equation}

The conditions (\ref{Contradiction_asumptions}) imply in particular the boundedness
\begin{equation}
\sup\limits_m~\Big(\| u^m\|_{L_3(Q^+)} + \| h^m \|_{L_3(Q^+)} + \| q^m \|_{L_{\frac 32}(Q^+)}\Big) \ < \ +\infty
\label{Basic_boundedness}
\end{equation}
From the local energy inequality near the boundary and the relation obtained from (\ref{MHD_Magnetic})
multiplied by the test function $\psi=\zeta (H^m)$ we obtain
\begin{equation}
\gathered
\| u^m \|_{L_{2,\infty}(Q^+_0(\frac 9{10}))} + \| h^m \|_{L_{2,\infty}(Q^+_0(\frac 9{10}))} \ + \\
+ \  \| u^m \|_{W^{1,0}_{2}(Q^+_0(\frac 9{10}))} + \| h^m \|_{W^{1,0}_{2}(Q^+_0(\frac 9{10}))} \ \le C.
\endgathered
\label{Gradient boundedness}
\end{equation}
From the equations (\ref{MHD_NS_m}), (\ref{MHD_Magnetic_m}) we also obtain the estimate
$$
\| \cd_t u^m \|_{L_{\frac 53}(-1, 0; W^{-1}_{\frac 53} (B^+_0)) } + \| \cd_t h^m \|_{L_{\frac 53}(-1, 0; W^{-1}_{\frac 53} (B^+_0)) }\ \le \ C.
$$

Hence we can extract subsequences
\begin{equation}
\begin{array}c
u^m \ \rightharpoonup \ u \quad \mbox{in}  \quad L_3(Q^+_0), \\
h^m \ \rightharpoonup \ h \quad \mbox{in}  \quad L_3(Q^+_0), \\
q^m \ \rightharpoonup \ q \quad \mbox{in}  \quad L_{\frac 32}(Q^+_0), \\
\end{array}
\label{Weak_Convergence0}
\end{equation}
\begin{equation}
\begin{array}c
u^m \ \rightharpoonup \ u \quad \mbox{in}  \quad W^{1,0}_2(Q^+_0(\frac 9{10})), \\
h^m \ \rightharpoonup \ h \quad \mbox{in}  \quad W^{1,0}_2(Q^+_0(\frac 9{10})),
\end{array}
\label{Weak_Convergence}
\end{equation}
\begin{equation}
\begin{array}c
u^m \ \to \ u \quad \mbox{in}  \quad L_3 (Q^+_0( \frac 9{10})), \\
h^m \ \to \ h \quad \mbox{in}  \quad L_3(Q^+_0( \frac 9{10})), \\
\vphi^m \to 0 \quad \mbox{in}  \quad C^2(\Bdva),\\
b_R(H^m)\ \to \  a  \quad \mbox{in}  \quad C^2,\qquad \  \qquad
\end{array}
\label{Strong_Convergence}
\end{equation}
here $a \in \R^3$ is the constant vector.

Passing to the limit in the equations (\ref{MHD_NS_m}), (\ref{MHD_Magnetic_m}) we obtain
\begin{equation}
\gathered
\cd_t u   - \Delta u +  \nabla q = \rot h \times   a \qquad \mbox{in}\quad Q^+_0,
\\
 \div u =0 \qquad \mbox{in}\quad Q^+_0,
 \\
u|_{y_3=0}=0,
\endgathered
\label{LinProb1}
\end{equation}
\begin{equation}
\gathered
\cd_t h     -  \Delta h = \rot (u\times  a) \qquad \mbox{in}\quad Q^+_0,
\\
\div h =0 \qquad \mbox{in}\quad Q^+_0,
\\
\begin{array}c
h_3|_{y_3=0}=0, \qquad \frac{\cd h_1}{\cd y_3}\big|_{y_3=0}\ =\ \frac{\cd h_2}{\cd y_3}\big|_{y_3=0} \ = \ 0.
\end{array}
\endgathered
\label{LinProb2}
\end{equation}
For the solution to the linear problem (\ref{LinProb1}) --- (\ref{LinProb2}) by a standard way (see \cite{VyaShi} Theorem 4.1) we obtain
\begin{equation}
Y_\tau(u) + \tilde Y_\tau(h) \ \le \ C(M)~ \tau^{1/3} ~Y_{1} (u, h, q)
\label{Y_u_h}
\end{equation}
Moreover from the second relation in \eqref{Contradiction_asumptions} we have
\begin{equation}
\liminf_{m \to \infty} Y_\tau (u^m,p^m,h^m) \geq c \tau^{\frac13}.
\label{lmc1}
\end{equation}
On the other hand we will show that
\begin{equation}
\limsup_{m \to \infty} Y_\tau (u^m,p^m,h^m) \leq c_* \tau^{\frac13}
\label{lmc2}
\end{equation}
for some constante $c_*$ Taking in \eqref{lmc1} a constant $c > c_*$ we obtain a contradiction. This
contradiction will prove the theorem.

From  (\ref{Strong_Convergence}) we conclude
$$
\lim\limits_{m\to +\infty} Y_\tau(u^m) = Y_\tau(u), \qquad \lim\limits_{m\to +\infty} \tilde Y_\tau(h^m) = \tilde Y_\tau(h)
$$
and hence
\begin{equation}
\limsup\limits_{m\to \infty}  Y_\tau (u^m, h^m , q^m ) \ \le \ Y_\tau(u) + \tilde Y_\tau(h) + \limsup\limits_{m\to \infty}  \hat Y_\tau (q^m).
\label{Y estimate}
\end{equation}
Then to prove \eqref{lmc2} it is sufficient to show that
\begin{equation}
\limsup_{m \to \infty} \hat Y_\tau (q^m) \leq c(M) \tau^{\frac13}.
\label{lcm3}
\end{equation}
For this purpose we decompose $(u^m, q^m)$ and $(u,q)$ as
$$
\gathered
u^m = u^m_1 + u^m_2, \qquad q^m = q^m_1+q^m_2, \\
u=u_1+u_2,\ \qquad\quad q = q_1+q_2,
\endgathered
$$
where $(u^m_1, q^m_1)\in W^{2,1}_{\frac 98, \frac 32}(\Q_0)\times W^{1,0}_{\frac 98, \frac 32}(\Q_0)$
are determined as a solutions of the following initial boundary-value problems in $\Q_0$:
$$
\gathered
\cd_t u^m_1 - \tilde \Delta_{\vphi_m} u^m_1 + \tilde\nabla_{\vphi_m} q_1^m = f^m \qquad \mbox{in}\quad \Q_0, \\
\tilde\nabla_{\vphi_m} \cdot u^m_1 =0 \qquad \mbox{in}\quad \Q_0, \\
u^m_1|_{t=-1} =0, \qquad u^m_1|_{y_3=0}=0,
\endgathered
$$
where $f^m$ is defined by the expression $\L_m \rot \L_m^{-T} h^m \times ( \delta_m h^m +  b_1(H^m)) -\delta_m (u^m\cdot \tilde\nabla_{\vphi^m} )u^m$
on the set $\Q_0(\frac 9{10})$
and extended by zero onto the whole $\Q_0$. Similarly,   $(u_1, q_1)$ are determined by the relations
\begin{equation}
\gathered
\cd_t u_1 - \Delta u_1 + \nabla q_1 = f \qquad \mbox{in}\quad \Q_0, \\
\div u_1 =0 \qquad \mbox{in}\quad \Q_0, \\
u_1|_{t=-1} =0, \qquad u_1|_{y_3=0}=0,
\endgathered
\label{Global_Linear_Problem}
\end{equation}
where $f$ determined by the expression $\rot h\times a $ on the set $\Q_0(\frac 9{10})$ and extended by zero onto the whole $\Q_0$.

As functions $u^m_1-u_1$, $q^m_1-q_1$ are the solution of the first initial boundary-value problem in $\Q_0$
with the right-hand side $f^m-f$ and zero initial and boundary conditions from lemma \ref{nonzerodivStokes_Lemma}, we obtain the estimate

\begin{equation}
\gathered
\| u^m_1 \|_{W^{2,1}_{\frac 98, \frac 32}(\Q_0)} +
\| \nabla q^m_1  \|_{L_{\frac 98, \frac 32}(\Q_0) } \ \le C \| f^m\|_{L_{\frac 98, \frac 32}(\Q_0(\frac 9{10}))}
\\
\| u^m_1-u_1 \|_{W^{2,1}_{\frac 98, \frac 32}(\Q_0)} + \| \nabla q^m_1 - \nabla q_1 \|_{L_{\frac 98, \frac 32}(\Q_0) } \
\le C \| f^m-f\|_{L_{\frac 98, \frac 32}(\Q_0(\frac 9{10}))}
\endgathered
\label{u_1 convergence}
\end{equation}
Note that
\begin{equation}
\gathered
\| f^m\|_{L_{\frac 98, \frac 32}(\Q_0(\frac 9{10}))} \ \le \ C(M)
\\
\| f^m - f\|_{L_{\frac 98, \frac 32}(\Q_0(\frac 9{10}))} \ \to \ 0, \quad \mbox{as} \quad m\to \infty.
\endgathered
\label{f^m boundedness}
\end{equation}
So, taking into account the imbedding $W^{1,0}_{\frac 98, \frac 32}(\Q_0(\frac 9{10})) \hookrightarrow L_{ \frac 32}(\Q_0(\frac 9{10}))$  we can conclude that
$$
\begin{array}c
q^m_1\to q_1\quad \mbox{in} \quad L_{\frac 32}(\Q_0(\frac 9{10}))
\end{array}
$$
and hence for any $\tau \in (0, \frac 9{10})$
$$
\lim\limits_{m\to \infty } Y_\tau (q^m_1 ) \ =  \ Y_\tau (q_1).
$$
On the other hand, $(u_1, q_1)$ is a solution of the linear Stokes problem in $\Q_0$. Hence from lemma \ref{zerodivStokes_Lemma} we conclude
$$
Y_\tau (q_1) \ \le \  C(M)~\tau^{1/3}~ Y_{\frac 9{10}} (q_1)
$$

We need to estimate $Y_{\frac 9{10}} (q_1)$.
From imbedding theorem $L_{\frac 32}(B^+_0(\frac{9}{10})) \hookrightarrow W^1_{\frac 98}(B^+_0(\frac{9}{10}))$ we conclude
$$
Y_{\frac 9{10}} (q_1)\ \le \ C ~\| \nabla q_1\|_{L_{\frac 98,\frac 32}(B^+_0(\frac{9}{10}))}
$$
For the solution $(u_1, q_1)$ of the initial-boundary value problem (\ref{Global_Linear_Problem}) we have the estimate
$$
\| u_1\|_{W^{2,1}_{\frac 98, \frac 32}(\Q_0(\frac 9{10}))} + \| \nabla q_1\|_{L_{\frac 98, \frac 32}(\Q_0(\frac 9{10}))} \
 \le \ C(M)~ \| \nabla h\|_{L_{\frac 98, \frac 32}(\Q_0(\frac 9{10}))}
$$
Using H\" older inequality \ $\| \nabla h\|_{L_{\frac 98, \frac 32}(\Q_0(\frac 9{10}))} \ \le \ C~\| \nabla h\|_{L_2(\Q_0(\frac 9{10}))}$
and  taking into account the weak convergence (\ref{Weak_Convergence}) from which we conclude
$$
\| \nabla h\|_{L_2(\Q_0(\frac 9{10}))} \ \le \ \liminf\limits_{m\to \infty} \| \nabla h^m \|_{L_2(\Q_0(\frac 9{10}))},
$$
and using (\ref{Gradient boundedness}) we obtain
$$
Y_{\frac 9{10}} (q_1) \ \le \ C(M).
$$

Now we consider functions $(u^m_2, q^m_2)$  determined by relations
\begin{equation}
\gathered
u^m_2:=u^m-u^m_1, \qquad q^m_2:= q^m -q^m_1.
\endgathered
\label{u_2 definition}
\end{equation}
These functions satisfy the homogeneous Stokes problems in $\Q_0(\frac{9}{10})$:
$$
\begin{array}c
\cd_t u^m_2 - \Delta u^m_2 + \nabla q_2^m = 0 \qquad \mbox{in}\quad \Q_0(\frac{9}{10}), \\
\div u^m_2 =0 \qquad \mbox{in}\quad \Q_0(\frac{9}{10}), \\
u^m_2|_{x_3=0}=0,
\end{array}
$$
$$
\begin{array}c
\cd_t u_2 - \Delta u_2 + \nabla q_2 = 0 \qquad \mbox{in}\quad \Q_0(\frac{9}{10}), \\
\div u_2 =0 \qquad \mbox{in}\quad \Q_0(\frac{9}{10}), \\
u_2|_{x_3=0}=0.
\end{array}
$$
Then
$$
\begin{array}c
\| u^m_2\|_{W^{2,1}_{9, \frac 32}(\Q_0(\frac 4{5}))}+ \| \nabla q^m_2\|_{L_{9, \frac 32}(\Q_0(\frac 4{5}))} \ \le \ C ~\Big(
\| u^m_2\|_{L_3(\Q_0(\frac 9{10}))} + \| q^m_2\|_{L_{\frac 32}(\Q_0(\frac 9{10}))}\Big)
\end{array}
$$
Note that due to (\ref{u_2 definition}), (\ref{Basic_boundedness}) and the first inequalities in (\ref{u_1 convergence}), (\ref{f^m boundedness})
we have the estimate
$$
\gathered
\| u^m_2\|_{L_3(\Q_0(\frac 9{10}))} + \| q^m_2\|_{L_{\frac 32}(\Q_0(\frac 9{10}))} \ \le \\
\le \
\| u^m\|_{L_3(\Q_0(\frac 9{10}))} + \| q^m\|_{L_{\frac 32}(\Q_0(\frac 9{10}))} \ +  \
\| u^m_1\|_{L_3(\Q_0(\frac 9{10}))} + \| q^m_1\|_{L_{\frac 32}(\Q_0(\frac 9{10}))}  \ \le  \\ \le \ C(M)
\endgathered
$$
On the other hand, from the H\" older inequality and lemma \ref{zerodivStokes_Lemma} we obtain for any $\tau \in (0, \frac 45)$
$$
\gathered
\hat Y_{\tau}(q^m_2)= \tau \Big( \frac1{|\Q_0(\tau)|} \int\limits_{\Q_0(\tau)} |q_2^m- [q_2^m]_{\Q_0(\tau)}|^{\frac 32}~dx dt\Big)^{\frac 23} \ \le \\
C \tau^2 ~ \Big( \frac1{|\Q_0(\tau)|} \int\limits_{\Q_0(\tau)} |\nabla q_2^m |^{\frac 32}~dx dt\Big)^{\frac 23} \\
\le \ C~\tau^{\frac 76} ~\| \nabla q^m_2 \|_{L_{9, \frac 32}(\Q_0(\frac 45))} \ \le \ C(M)~\tau^{\frac 76}
\endgathered
$$

Summarizing all previous estimates we arrive at
$$
\gathered
\limsup\limits_{m\to \infty}  \hat Y_\tau (q^m) \ \le \  \lim\limits_{m\to \infty}  \hat Y_\tau (q^m_1) + \limsup\limits_{m\to \infty}  \hat Y_\tau (q^m_2) \
\le \ C(M) ~ \tau^{\frac 13}
\endgathered
$$
which gives us a contradiction with \eqref{lmc1}.

\qed

Iterating \eqref{Decay_Estimate} and using scaling argument it is easy to obtain the following lemma (see \cite{ESS} and \cite{SSS} for details).

\begin{lemma}
There exists an absolute constant $\ep_{**}>0$ such that for any $M>0$ and $\beta \in (0,1/3)$ there exists $\tau \in (0,1/2)$ with the following properties.
For any boundary suitable weak solution $(v,H,p,\vphi)$ of the MHD system (\ref{MHD_NSE}), (\ref{MHD_Magnetic}) near the boundary in $\Q(1)$
the following implication holds:
if
$$
Y_1(v,H,p) + \| \vphi \|_{C^2} \ < \ \ep_{**},
$$
and
$$
b_R(H) \le \  M
$$
then
\begin{equation}
Y_{\tau^k} (v,H,p) \ \le \tau^{\beta k}~Y_1 (v,H,p)
\end{equation}
\end{lemma}

Theorem \ref{CKN_theorem} follows from this lemma in the standard way
by scaling arguments, and combination of boundary estimates with the internal estimates obtained in \cite{Vya}.
See  details in  \cite{Seregin_JMFM}, \cite{Seregin_Aa}, \cite{SSS}, \cite{Seregin_Handbook}.

\section{Estimates of solutions of the heat equation with homogeneous boundary data}
\label{Heat_equation_estimate_section}
\setcounter{equation}{0}

In this section we will obtain some estimates for $L_2$-norms of solutions of homogeneous initial and boundary problem for the heat equation in half-ball.
Namely, we consider the following problem
\begin{equation}
\gathered
\cd_t h - \lap h = f \quad \text{ in } \Q_0(R),\\
h_3|_{x_3=0} = 0, \quad \frac{\cd h_i}{\cd x_3}|_{x_3=0} = 0 \quad \text{in } \Qdva(R) \quad i=1,2,\\
h|_{t=-R^2} = 0,
\endgathered
\label{rhsel1}
\end{equation}
here $h: \Q_0(R) \to \R^3$ is an unknown function.

The main result of this section is the following theorem
\begin{theorem}
Let $f \in L_{\frac 32,1}( \Q_0(R) )$, and $h$ is the solution of \eqref{rhsel1}. Then the following estimate holds
\begin{equation}
\| h \|_{L_2( \Q_0(R) )} \leq c R^{\frac 12} \| f \|_{L_{\frac 32, 1}}( \Q_0(R) ).
\label{lh1}
\end{equation}
\label{lhTh1}
\end{theorem}

We note, that conditions for $h$ on a plain part of a boundary allow us to extend this function
into whole $B$ by the following way: components $h_1$ and $h_2$ will be extended as even functions and component $h_3$ as odd function. The right hand side
can be extended by the same manner. We also put $f \equiv 0$ in $\R^3 \backslash \B(R)$.
So it is sufficient to prove the theorem for the solution of the following Cauchy problem for the heat equation.
\begin{equation}
\gathered
\d_t h - \Delta{h} = f \quad \text{ in } \quad \Pi_R, \\
h|_{t= -R^2} = 0,
\endgathered
\label{Heat_eq}
\end{equation}
here $\Pi_R = \R^3 \times [-R^2,0)$.

To prove this theorem we will need the Young inequality for convolutions (se. \cite{Hardy}, \cite{Stein}). Namely, let
$$ g(x) = \intl_{\R^n} K(x-y) f(y)\, dy,$$
then for arbitrary $1 \leq p \leq q \leq \infty$, the following estimate holds
\begin{equation}
\| g \|_q \leq \| K \|_l \| f \|_p, \quad \text{ here } \quad 1 - \frac1p + \frac 1q = \frac 1l.
\label{convYung}
\end{equation}
In particular we will use an inequality
\begin{equation}
\| g \|_2 \leq \| K \|_{\frac65} \| f \|_{\frac 32}.
\label{lh3}
\end{equation}

\begin{lemma}
Let $f \in L_1(-R^2,0)$, for some $R > 0$ and
\begin{equation}
g(t) = \intl_{-R^2}^t \frac{ f(\tau)\, d\tau}{ ( t - \tau )^{\frac14}},
\label{lh3}
\end{equation}
then
\begin{equation}
\| g \|_{L_2(-R^2,0)} \leq c R^{\frac12} \| f \|_{L_1(-R^2,0)}
\label{lh4}
\end{equation}
\end{lemma}

{\bf Proof.} In the case $R = 1$ inequality \eqref{lh4} is the corollary of \eqref{convYung} with $p = 1$ and $q = l = 2$.
In general case we make scaling transformations.
Namely, let consider the functions
$$f^*(s) = f (R^2 s), \quad g^*(s) = g(R^2 s), \quad \hat{g}(s) = \intl_{-R^2}^s \frac{ f^*(\sigma)\, d\sigma}{(s - \sigma)^{\frac14}}.$$
Changing variables under the integral, we have
$$
\gathered
g^*(s) = R^{\frac32} \hat{g}(s), \\
\|f^* \|_{L_1(-1,0)} = R^{-2} \|f \|_{L_1(-R^2,0)}, \quad \| g \|_{L_2(-R^2,0)} = R \| g^* \|_{L_2(-1,0)}.
\endgathered
$$
And then we get
$$
\gathered
\| g \|_{L_2(-R^2,0)} = R \| g^* \|_{L_2(-1,0)} = R^{\frac52} \| \hat{g} \|_{L_2(-1,0)} \leq \\
\leq c R^{\frac52} \| f^* \|_{L_1(-1,0)} = c R^{\frac12} \| f \|_{L_1(-R^2,0)}.
\endgathered
$$

\qed

{\bf Proof of theorem \ref{lhTh1}.} The solution of \eqref{Heat_eq} can be found as follows
$$ h(t,x) = \intl_{-R^2}^t \intl_{\R^3} \frac{ e^{-\frac{|x-y|^2}{4(t-\tau)}}}{(4 \pi (t - \tau))^{\frac32}} f(\tau,y) \,dy d\tau.$$
We fix $t$, look on excretion under the integral by time as function with values in Banach space $L_2(\R^3)$ and use
inequality to the norm of its integral
\begin{equation}
\| h (t,\cdot) \|_{L_2(\R^3)} \leq
 \intl_{-R^2}^t \l\| \intl_{\R^3} \frac{ e^{-\frac{|x-y|^2}{4(t-\tau)}}}{(4 \pi (t - \tau))^{\frac32}} f(\tau,y) \,dy\r\|_{L_2(\R^3)} \,d\tau.
\label{lh5}
\end{equation}
By direct computations we find
\begin{equation}
\l( \intl_{\R^3} e^{-\frac{z^2}{4(t-\tau)} \cdot \frac65} \,dz \r)^{\frac56} = c (t - \tau )^{\frac54}.
\label{lh6}
\end{equation}
Then from \eqref{lh3}, \eqref{lh5} and \eqref{lh6} we obtain
$$\| h(t,\cdot) \|_{L_2(\R^3)} \leq c \intl_0^t \frac{\| f(\tau,\cdot) \|_{L_{\frac32}(\R^3)}}{(t - \tau)^{\frac14}} \,d\tau.$$
Next form \eqref{lh4}
$$ \|h \|_{L_2(\Pi_R)} \leq c R^{\frac12} \| f \|_{L_{\frac32 , 1}(\Pi_R)}.$$

\qed


\section{Estimates for the magnetic component}
\label{H_estimate_section}
\setcounter{equation}{0}

In this section we will get estimates for magnetic component of the suitable weak solution in $\Q(R)$.

We recall, that $H$ satisfy the following integral identity
\begin{equation}
\intl_{\Q(R)} ~ \Big( - H\cdot \cd_t \psi  + \rot H\cdot \rot \psi - (v\times H) \cdot \rot \psi \Big)~dxdt  \ = 0,
\label{lhes1}
\end{equation}
for all \ $\psi \in C^{\infty}(\Q(R))$ \ such that \  $\psi_\nu |_{\cd\Om\times (-R^2,0)}=0$ and $\psi(-R^2,x) = \psi(0,x) = 0$.
Also we have $H_\nu |_{\Gamma_R \times (-R^2,0)}=0$ and $\div H = 0$.
Note, that without loss of generality we can assume, that \eqref{lhes1} holds only for test functions with $\div \psi = 0$.

\begin{theorem}
Assume that \eqref{lhes1} holds for some function $v \in \WO$.
Then there exist absolute positive constants $\ep_1$, $\al$ and $c$ such that for any $\ep\in (0,\ep_1)$ and any $K>0$ if
\begin{equation}
\sup\limits_{r\in (0,1)} E(r) < \ep, \quad \| \vphi \|_{C^2} < \ep \quad \text{and} \quad \sup\limits_{r\in (0,1)} E_*(r) < K
\label{Assumptions}
\end{equation}
then for any $0<r< R \leq 1$
\begin{equation}
 F_2(r) \ \le \   c\left(\frac{r}{R}\right)^{2} F_2(R) + c \ep (F_2(R) + K + 1).
\label{F_2_l1}
\end{equation}
\label{Bound_F_2}
\end{theorem}

{\bf Proof.} We proceed to the coordinates \eqref{change_variables}.
Then \eqref{lhes1} transforms to
$$
\intl_{-R^2}^0 \l[ -(\tH,\d_t\tpsi) + ( \L \rot \L^T \tH , \L \rot \L^T \tpsi) - ( \tv \times \tH , \L \rot \L^T  \tpsi) \r]\,dt = 0,
$$
here $(\cdot , \cdot)$ is $L_2$ inner product.

Next we introduce new functions $h = \L^{-1} \tH$ and $\eta = \L^{-1} \tpsi$. Then the last identity can be written as follows
\begin{equation}
\intl_{-R^2}^0 \l[ -(h,\d_t \A \eta) + (\A \rot \A h,\rot \A \eta) - (\tv \times \tH, \L \rot \A \eta) \r]\,dt = 0
\label{lhes2}
\end{equation}
here $\A = \L^T \L$.

Note, that $\nabla_y \cdot h = \nabla_y \cdot \L^{-1} \tH = \L^{-T} \nabla_y \cdot \tH = \nabla_x \cdot H \circ e = 0$.
Also $H \cdot \nu \circ e = \L h \cdot \tilde{\nu} = h \cdot \L^{T} \tilde{\nu} = -h_3$. Similar identities holds for function $\eta$.

As the result we can consider function $h$ as the generalized solution of parabolic system which corresponds to holding the identity \eqref{lhes2}
for arbitrary function $\eta \in C^{\infty}$, such that
\begin{equation}
\div \eta = 0 \quad \text{ and } \quad \eta_3|_{y_3 = 0} = 0.
\label{lhes5}
\end{equation}

Then to estimate $L_2$-norm of $H$ it will sufficient to obtain
inequality for $\| h \|_2$. To do this we decompose it into three
parts
\begin{equation}
h = \hi + \hii + \hiii.
\label{lhes7}
\end{equation}

Here $\hi$ is the solution of the following initial-boundary problem
\begin{equation}
\begin{aligned}
&\d_t \hi - \lap \hi = \A \rot \l( \L^T ( \tv \times \tH) \r)\\
& \hi_3|_{y_3 = 0} = 0, \quad \frac{\d \hi_1}{\d y_3} |_{y_3 = 0} = 0, \quad \frac{\d \hi_2}{\d y_3} |_{y_3 = 0} = 0,\\
&\hi|_{t = -R^2} = 0.
\end{aligned}
\label{lhes3}
\end{equation}
From \eqref{rhsel1} we have
$$ \| \hi \|_{2,R} \leq c R^{\frac12} \| \A \rot \l( \L^T ( \tv \times \tH) \r) \|_{\frac32,1}
\leq c(\vphi) R^{\frac12} ( \| |v||\nabla H| + |\nabla v||H| \|_{\frac32,1}).$$
Next we use H\"older inequality and the embedding theorem
\begin{equation}
F_2(R,\hi) \leq c(\vphi) \l( E(R)E_*(R) + E(R)F_2(R) \r).
\label{lhes4}
\end{equation}
Boundary conditions for function $\hi$ imply that the following identity
\begin{equation}
\intl_{-R^2}^0 \l[ -(\hi,\d_t \eta) + (\rot \hi,\rot \eta) - (\tv \times \tH, \L \rot \A \eta) \r]\,dt = 0
\label{lhes6}
\end{equation}
for every function $\eta$ satisfying \eqref{lhes5}.
As \eqref{lhes4} and \eqref{lhes6} are stored at the replacing $\hi$ to the solenoidal component of its Weil decomposition
without loss of generality we can assume, that $\div \hi = 0$.

The second component of \eqref{lhes7} is the solution of the following problem
\begin{equation}
\gathered
\intl_{-R^2}^0 \l[ -(\hii,\d_t \eta) + (\rot \hii, \rot \eta)\r]\,dt =\\
\intl_{-R^2}^0 \l[ (h, \d_t( \A - I)) + (\rot h, \rot \eta) - (\A \rot \A h, \rot \A \eta) \r]\,dt
\endgathered
\label{lhes8}
\end{equation}
for arbitrary function $\eta$ satisfying \eqref{lhes5} with the initial and boundary conditions
$$
\hii_3|_{y_3 = 0} = 0, \quad \hii|_{t = -R^2} = 0.
$$
To obtain the estimate for $\hii$ we consider the dual problem
\begin{equation}
\begin{aligned}
&  \d_t \eta + \lap \eta = -\hii \\
& \eta_3|_{y_3 = 0} = 0, \quad \frac{\d \eta_1}{\d y_3} |_{y_3 = 0} = 0, \quad \frac{\d \eta_2}{\d y_3} |_{y_3 = 0} = 0,\\
&\eta|_{t = 0} = 0.
\end{aligned}
\label{lhes9}
\end{equation}
For the right hand side of \eqref{lhes8} we have the following identity
\begin{equation}
\gathered
(\rot h , \rot \eta) - (\A \rot \A h, \rot \A \eta) = \\
( \rot h , \rot (I - \A) \eta) + (\rot (I - \A)h, \rot \A \eta) + ((I - \A) \rot \A h, \rot \A \eta).
\endgathered
\label{lhes10}
\end{equation}
Substituting into \eqref{lhes8} the solution of \eqref{lhes9} we have
$$
\| \hii \|_2 = \intl_{-R^2}^0 \l[ (h, \d_t( \A - I)) + (\rot h, \rot \eta) - (\A \rot \A h, \rot \A \eta) \r]\,dt
$$
Note that the matrix $\A$ is close the identity, so from \eqref{lhes10} and coercive estimates for \eqref{lhes8} we obtain
\begin{equation}
F_2(R, \hii) \leq c \| \vphi \|_{C^2} \l( F_2(R) + E_*(R) \r).
\label{lhes11}
\end{equation}

The third component of \eqref{lhes7} satisfy to the homogenous boundary problem for the heat equation
\begin{equation}
\begin{aligned}
&  \d_t \hiii - \lap \hiii = 0 \\
& \hiii_3|_{y_3 = 0} = 0, \quad \frac{\d \hiii_1}{\d y_3} |_{y_3 = 0} = 0, \quad \frac{\d \hiii_2}{\d y_3} |_{y_3 = 0} = 0.
\end{aligned}
\label{lhes12}
\end{equation}
Extending $\hiii$ into whole cylinder and using mean value theorem we have
\begin{equation}
\gathered
F_2(r,\hiii) \leq c \l(\frac{r}{R} \r)^2 F_2 (R, \hiii) \\
\leq c \l(\frac{r}{R} \r)^2 \l( F_2 (R, h) + F_2(R,\hi + \hii) \r).
\endgathered
\label{lhes13}
\end{equation}

Combining \eqref{lhes4}, \eqref{lhes11} and \eqref{lhes13} we obtain the statement of the theorem.

\qed


\section{Estimates of Energy Functionals}
\label{EnergyEstimates_Section}
\setcounter{equation}{0}

Now we define few more functionals.
Note that all these functionals are invariant with respect to the natural scaling of the MHD system.
For $r\le 1$, $q\in [1,\frac{10}3]$, $s\in [1,\frac 98]$ and $(v,p,H,\vphi)$ suitable weak solution to the MHD system in $\Q(R)$ $0 < r < R <1$
we introduce the following quantities:
$$
\begin{array}c
A( r) \equiv \Big( \frac 1{r}
\sup\limits_{t\in (-r^2, 0)} \int\limits_{B^+(r)} |v|^{2}~dy
\Big)^{1/2}, \\
A_*( r) \equiv \Big( \frac 1{r}
\sup\limits_{t\in (-r^2, 0)} \int\limits_{B^+(r)} |H |^{2}~dy
\Big)^{1/2},
\\
C_q( r) \equiv \Big( \frac 1{r^{5-q}} \int\limits_{Q^+(r)} |v|^q~dydt
\Big)^{1/q},
\\ D(r) \equiv \Big( \frac 1{r^2} \int\limits_{Q^+(r)}
|p - [p]_{B^+(r)}|^{3/2}~dydt \Big)^{2/3},
\\
D_s(r) = R^{\frac 53 -\frac 3s} \Big( \int\limits_{-r^2}^{0} \Big(
\int\limits_{B^+(r)}  |\nabla p|^{s}~dy \Big)^{\frac 1s \cdot
\frac 32}~dt \Big)^{2/3},
\end{array}
$$
$$
C(r) = C_3(r), \qquad F(r)= F_3(r), \qquad D_*(r)= D_{\frac{36}{35}}(r).
$$

\noindent
First we formulate the set of results following from the general theory of functions:

\begin{lemma}
\label{Interpolation}
Let $R > 0$, $\vphi \in C^2(\Bdva(R))$, $v$, $H\in W^{1,0}_2(Q^+(R))$ and $p\in W^{1,0}_{\frac 98, \frac 32}(Q^+(R))$ are arbitrary functions.
Assume $v|_{x_3= \vphi(x_1,x_2)}=0$. Then for any $0 < r <R$  the following inequalities hold:
\begin{equation}\label{C_3}
C(r) \ \le  \ ~A^{\frac 12}(r)E^{\frac 12}(r),
\qquad
F(r)\ \le \ A_*^{\frac 12}(r)[ E_*^{\frac 12}(r) +  F_2^{\frac 12}(r)]
\end{equation}
\begin{equation}
D(r)\ \le \  c D_1(r), \qquad D_1(r)\ \le \ c D_s(r), \qquad \forall s>1.
\label{D}
\end{equation}
\end{lemma}

The proof of this lemma follows from interpolation inequalities and imbedding theorems.
Proof of the similar inequalities for the Navier-Stokes system can be found in \cite{LS}.

\begin{lemma}
\label{EnergyEstimate_Lemma}
Assume $(v,p,H,\vphi)$ is a boundary suitable weak solution to the MHD equations in $Q^+$.
Then for any $r\in (0,1)$ the following inequality holds
\begin{equation}
\gathered
A(r/2) + A_*(r/2) + E(r/2) + E_*(r/2) \ \le \\ \le \ c~ \Big( C_2(r)+ F_2(r) + C^{\frac 12}(r)D^{\frac 12}(r) +  C^{\frac 32}(r)\Big) \ + \\
+ \ c~ \Big( C^{\frac 12}(r)A_*^{\frac 12}(r)E_*^{\frac 12}(r) + F^{\frac 12}(r)A_*^{\frac 12}(r)E^{\frac 12}(r) \Big)
\label{Lo_En_Iq}
\endgathered
\end{equation}
\end{lemma}

\noindent
{\bf Proof.} Estimate (\ref{Lo_En_Iq}) follows from (\ref{LEI}) in a standard way. We just explain the specific estimates of the terms
$$
I_1:=\int\limits_{Q^+(r)} ~ |H|^2 (v \cdot \nabla \zeta)   ~dxdt  \quad \mbox{and}   \quad I_2:=\int\limits_{Q^+(r)} ~ (v \cdot H) (H\cdot \nabla \zeta)  ~dxdt.
$$

$I_1$ we transform in the following way
$$
\begin{array}c
I_1\ = \ \int\limits_{Q^+(r)} ~ \Big( |H|^2  - [|H|^2]_{B^+(r)} \Big) (v \cdot \nabla \zeta)  ~dxdt
\end{array}
$$
Applying the H\" older  inequality we obtain
$$
|I_1| \ \le \ \frac cr ~\int\limits_{-r^2}^0 \left\| |H|^2  - [|H|^2]_{B^+(r)}\right\|_{L_{\frac 32}(B^+(r))} \| v\|_{L_{3}(B^+(r))}~dt
$$
Applying the inequality $\| f-[f]_{B^+(r)}\|_{L_{\frac 32}(B^+(r))} \le c \| \nabla f \|_{L_1(B^+(r))}$, we arrive at
$$
\gathered
|I_1| \ \le \ \frac cr ~\int\limits_{-r^2}^0 \| \nabla |H|^2  \|_{L_{1}(B^+(r))} \| v\|_{L_{3}(B^+(r))}~dt  \ \le \\ \le \
\frac cr ~\int\limits_{-r^2}^0 \| H\|_{L_2(B^+(r))} \| \nabla H  \|_{L_{2}(B^+(r))} \| v\|_{L_{3}(B^+(r))}~dt  \ \le \\ \le \
\frac c{r^{2/3}}  ~\| H\|_{L_{2,\infty}(Q^+(r))} \| \nabla H \|_{L_2(Q^+(R))} \| v \|_{L_3(Q^+(r))}\ \le \  cr~ A_*(r) E_*(r) C(r)
\endgathered
$$

For $I_2$ we obtain relations
$$
I_2 \ = \ \int\limits_{Q^+(r)} ~ \Big ( (v \cdot H)- [v \cdot H]_{B^+(r)}\Big) (H\cdot \nabla \zeta)  ~dxdt
$$
Hence
$$
\gathered
|I_2| \ \le   \ \frac cr ~\int\limits_{-r^2}^0 \left\| (v \cdot H)- [v \cdot H]_{B^+(r)} \right\|_{L_{2}(B^+(r))} \| H\|_{L_{2}(B^+(r))}~dt \ \le \\
\le \ \frac cr ~ \| H\|_{L_{2,\infty}(Q^+(r))} \int\limits_{-r^2}^0 \left\| \nabla (v \cdot H) \right\|_{L_{\frac 65}(B^+(r))} ~dt \ \le \
\frac cr ~ \| H\|_{L_{2,\infty}(Q^+(r))} \ \times \\ \times \ \int\limits_{-r^2}^0 \Big( \| \nabla v \|_{L_2(B^+(r))} \| H \|_{L_3(B^+(r))}
 + \| \nabla H \|_{L_2(B^+(r))} \| v \|_{L_3(B^+(r))}\Big)~dt \ \le \\
\le \ \frac c{r^{2/3}}  ~\| H\|_{L_{2,\infty}(Q^+(r))}  \Big( \| \nabla v \|_{L_2(Q^+(r))} \| H \|_{L_3(Q^+(r))}
+ \| \nabla H \|_{L_2(Q^+(r))} \| v \|_{L_3(Q^+(r))}\Big) \endgathered
$$
So, we obtain
$$
|I_2| \ \le \  cr~ A_*(r) ~\Big(~E(r) F(r) + E_*(r) C(r)~\Big)
$$

\qed

\begin{lemma}
\label{PresureEstimate_Lemma}
Assume $(v,p,H,\vphi)$ is a boundary suitable weak solution to the MHD equations in $Q^+$ and $\| \vphi \|_{C^2(\Bdva)} < \mu < \frac{\mu_*}2$
where $\mu_*$ is the constant defined in lemma \ref{nonzerodivStokes_Lemma}.
Then for any $r\in (0,1)$ and $\theta \in (0,\frac14)$ the following inequality holds
\begin{equation}
\gathered
D_*(\theta r) \ \le \ c ~\theta^{\frac 43} ~\Big( D_*(r) + E(r)  \Big) \ + \\
+  \ c(\theta) ~\Big(  A^{\frac 23}(r) E^{\frac 43}(r) +  A^{\frac 56}_*(r)  F^{\frac 16}(r) E_*(r) \Big)
\endgathered
\label{D_*}
\end{equation}
\end{lemma}

\noindent
{\bf Proof.}
To obtain (\ref{D_*}) we apply the method developed in \cite{Seregin_JMFM}, \cite{Seregin_ZNS271}, see also \cite{SSS}.
Let $e(y)$ is the map defined by \eqref{change_variables}.
We fix $r\in (0,1]$ and $\theta\in (0,\frac 14)$ and without loss of generality we can assume, that
$e^{-1}(\Q(\theta r) \subset \Q_0(2\theta r) \subset \Q_0(r/2) \subset e^{-1}(\Q(r))$.
Then we decompose $v$ and $p$ as
$$
v\ = \ \hat v + \check v,\qquad p \ = \ \hat p +\check p,
$$
where $(\hat v, \hat p)$ is a solution of the perturbed Stokes initial boundary value problem in a half-space
$$
\gathered
\left\{\begin{array}c  \cd_t \hat v - \tilde\Delta_{\vphi} \hat v+\tilde\nabla_{\vphi} \hat p \ = \L \rot (\L^{-T} \tH) \times \tH-(\tv \cdot \tilde\nabla_{\vphi} ) \tv , \\
\tilde\nabla_{\vphi} \hat v =0 \end{array}\right. \qquad \mbox{in}\quad \Q_0(\frac{r}2), \\
\hat v|_{t=0}=0, \qquad \hat v|_{y_3=0}=0,
\endgathered
$$
and $(\check v, \check p)$ is a solution of the homogeneous perturbed Stokes system in $Q_0^+(\frac{r}2)$:
$$
\gathered
\left\{\begin{array}c  \cd_t \check v - \tilde\Delta_{\vphi} \check v +\tilde\nabla_{\vphi} \check p\ = \ 0, \\
\tilde\nabla_{\vphi} \check v =0 \end{array}\right. \qquad \mbox{in}\quad Q^+_0(\frac{r}2),  \\
\check v|_{y_3=0}=0.
\endgathered
$$

For $\nabla \hat p$ and $\nabla \check p$ from lemmas \ref{nonzerodivStokes_Lemma} and \ref{zerodivStokes_Lemma} we have the following estimates.
$$
\gathered
\| \nabla \hat p\|_{L_{\frac{36}{35},\frac 32}(Q^+_0(\frac{r}2))} \ + \ \frac 1r \| \nabla \hat v \|_{L_{\frac{36}{35},\frac 32}(Q^+_0(\frac{r}2))}   \ \le \\
\le \  c~\Big( ~\| H\times \rot H\|_{L_{\frac{36}{35}, \frac 32}(Q^+(r))} \ + \ \| (v\cdot \nabla)v \|_{L_{\frac{36}{35}, \frac 32}(Q^+(r))}~\Big),
\endgathered
$$
$$
\| \nabla \check p\|_{L_{\frac{36}{35},\frac 32}(Q^+_0(\theta r))} \
\le \ c~\theta^{\frac{31}{12}} ~\Big( ~\frac 1r \| \nabla \check v \|_{L_{\frac{36}{35},\frac 32}(Q_0^+(\frac{r}2))}
\ + \  \|  \nabla \check p \|_{L_{\frac{36}{35},\frac 32}(Q_0^+(\frac{r}2))}~\Big).
$$

From the H\" older inequality we obtain
$$
\gathered
\| H\times \rot H\|_{L_{\frac{36}{35}, \frac 32}(Q^+(r))} \
\le \ c ~r^{\frac 29} ~\| H \|_{L_{2,\infty}(Q^+(r))}^{\frac 56} \| \nabla H\|_{L_2(Q^+(r))} \|  H \|_{L_3(Q^+(r))}^{\frac 16}
\endgathered
$$
$$
\gathered
\| (v\cdot \nabla)v \|_{L_{\frac{36}{35}, \frac 32}(Q^+(r))} \ \le \ c~r^{\frac 14}~ \| (v\cdot \nabla )v \|_{L_{\frac{9}{8}, \frac 32}(Q^+(r))} \ \le  \\
\le \ c~r^{\frac 14}~
\| v\|_{L_{2,\infty}(Q^+(r))}^{\frac 23} \| \nabla v\|_{L_2(Q^+(r))}^{\frac 43}
\endgathered
$$
Representing $\check v = v-\hat v$, $\check p=p-\hat p$ and gathering all  above estimates for $\hat p$ and $\hat v$ we obtain
$$
\gathered
D_*( \theta r) \ \le \ c~\theta^{\frac 43} ~\Big( ~D_*(r) + E(r) + A^{\frac 23}(r) E^{\frac 43}(r) + A_*^{\frac 56}(r) E_*(r) F^{\frac 16}(r) ~\Big) \ +  \\
+ \ c(\theta)~ \Big( ~A^{\frac 23}(r) E^{\frac 43}(r) + A_*^{\frac 56}(r) E_*(r) F^{\frac 16}(r) ~\Big)
\endgathered
$$

\qed


\section{CKN condition and Partial Regularity of Solutions}
\label{Final_section}
\setcounter{equation}{0}

In this section we present the proofs of Theorems \ref{CKN_theorem} and \ref{Partial_Regularity}.
We start from proof of the modified version of \eqref{F_2_l1}.

\begin{lemma}
For any $K>0$ there exists a constants $c(K)>0$ and $\ep_2>0$ such that  for any $\ep\in (0,\ep_2]$ and
any boundary suitable weak solution  $(v,H,p,\vphi)$ of the MHD system in $Q^+$ if
\begin{equation}
\sup\limits_{r\in (0,1)}  E(r)\le  \ep, \quad \| \vphi \|_{C^2(\Bdva)} < \ep, \quad  \sup\limits_{r\in (0,1)} E_*(r)  \ \le \ K,
\label{Statement_of_Theorem_E_E1}
\end{equation}
then for some $\al > 0$ and any $0<r< R \leq 1$
\begin{equation}
 F_2(r) \ \le \   c\left(\frac{r}{R}\right)^{\al} F_2(R) + c(K) \ep_2.
\label{F_2}
\end{equation}
\label{ckn_lemma1}
\end{lemma}

\noindent
{\bf Proof.} We will use standard iteration technic. Let $R > 0$ and $\theta \in (0,1/2)$. We fix $\ep_1$ from theorem \ref{Bound_F_2}.
Then from \eqref{F_2_l1} we have
\begin{equation}
F_2(\theta R) \leq c_1 (\theta^2 + \ep_1) F_2(R) + c \ep_1 (K+1).
\label{cknl1}
\end{equation}
Next we choose $\theta$ and $\ep_2 < \ep_1$ such that
$$ c_1 \theta^2 \leq \frac14, \quad c_1 \ep_2 \leq \frac14.$$
Then from \eqref{cknl1} we obtain
$$ F_2(\theta R) \leq \frac12 F_2(R) + c \ep_2 (K+1).$$
Next we will iterate the last inequality
\begin{equation}
\gathered
F_2 (\theta^k R) \leq \frac12 F_2(\theta^{k-1} R) + c \ep_2 (K+1) \leq\\
\leq \frac14 F_2(\theta^{k-2} R) + \l( 1 + \frac12 \r) c \ep_2 (K+1)\leq\\
\leq \frac1{2^k} F_2(R) + c \ep_2 (K+1).
\endgathered
\label{cknl2}
\end{equation}
Finally we put $\al = \log_{\theta} \frac12$ and chose $k>0$ such that $ \theta^{k+1} R \leq r \leq \theta^k R$.
Then from \eqref{cknl2} we obtain
$$
\gathered
F_2(r) \leq c F_2(\theta^k R) \leq c \frac{1}{2^k} F_2(R) + c(K) \ep_2 \leq\\
\leq c \theta^{k\al} F_2(R) + c(K) \ep_2 \leq c \theta^{-\al} \l( \frac{\theta^{k+1} R}{R} \r)^{\al} F_2(R) + c(K) \ep_2 \leq\\
\leq c \l( \frac{r}{R} \r)^{\al} F_2(R) + c(K) \ep_2.
\endgathered
$$

\qed


\begin{lemma} \label{Boundedness_A_A_D}
Denote by \ $\mathcal E(r)$ \ the following functional
$$
\mathcal E (r) \ = \ A(r)+   A_*(r) +  D_*(r),
$$
and let $\ep_2>0$ be the absolute  constant defined in lemma \ref{ckn_lemma1}.
For any $K>0$ there exists a constant $c(K)>0$ such that  for any $\ep\in (0,\ep_2]$ and
any boundary suitable weak solution  $(v,H,p,\vphi)$ of the MHD system in $Q^+$ if
\begin{equation}
\sup\limits_{r\in (0,1)}  E(r)\le  \ep, \quad \| \vphi \|_{C^2(\Bdva)} < \ep, \quad  \sup\limits_{r\in (0,1)} E_*(r)  \ \le \ K,
\label{Statement_of_Theorem_E_E1}
\end{equation}
and
\begin{equation}
F_2(1) \ \le \ M,
\label{M}
\end{equation}
then for any $0< r < R \le 1$
\begin{equation}
\mathcal E(r) \ \le \ c \left(\frac{r}{R}\right)^\be\mathcal E(R) \  +  \ c(K)(1+R^\al M).
\label{Estimate_varE}
\end{equation}
where $\be>0$ is some absolute constant.
\end{lemma}

\noindent
{\bf Proof.} Without loss of generality we can assume $K\ge 1$. Then from  (\ref{F_2})  we obtain
$$
F_2(R)\le cr^\al M  + c \ep_2  C(K).
$$
 From this inequality and (\ref{C_3}) we obtain
\begin{equation}
C(R) \ \le \ c~ \mathcal E^{\frac 12} (R) \ep_1^{\frac 12}, \qquad F(R)\
\le \ c~ \mathcal E^{\frac 12} (R)\Big(C(K) +  R^{\frac \al 2}M^{\frac 12} \Big)
\label{C}
\end{equation}

Assume $r\in (0,1)$ and $\theta \in (0,\frac 12)$.
From (\ref{Lo_En_Iq}) with the help of (\ref{D}) and the Young inequality  we obtain
$$
\gathered
\mathcal E(\theta R) \ \le \  c~\Big( F_2(2\theta R) + D_*(2\theta R)\Big) \  +
\\ + \
c(\theta)\Big( C_2(R) + C(R) + C^{\frac 32}(R) + C^{\frac 12}(R)A_*^{\frac 12}(R)E_*^{\frac 12}(R) + F^{\frac 12}(R)A_*^{\frac 12}(R)E^{\frac 12}(R)\Big)
\endgathered
$$
Taking into account (\ref{C}) and (\ref{Statement_of_Theorem_E_E1}) we obtain
\begin{equation}
\gathered
\mathcal E(\theta R) \ \le \  c~\Big( F_2(2\theta R) + D_*(2\theta R)\Big) \  +
\\ + \
c(\theta)\l(\mathcal E^{\frac 12} (R) \ep_1^{\frac 12} + \mathcal E^{\frac 34} (R) \ep_2^{\frac 34} + \ep_2^{\frac 14}\mathcal E^{\frac 34}(R)K^{\frac 12} + \r.\\
 \l.+ ( C(K) +  R^{\frac \al 4}M^{\frac 14})
 \mathcal E^{\frac 34}(R)\ep_1^{\frac 12}\r)
\endgathered
\label{Similar to this}
\end{equation}
Applying the Young inequality $ab\le \ep a^p+C_\ep b^{p'}$ we obtain
$$
\gathered
\mathcal E(\theta R)  \le   \frac 14 \mathcal E(R)  +   c\Big( F_2(2\theta R) + D_*(2\theta R)\Big)   +  c(\theta )c(K)  +   c(\theta) R^\al M.
\endgathered
$$

From (\ref{F_2}) and (\ref{D_*}) we obtain
$$
\gathered
F_2(2\theta R) + D_*(2\theta R)  \le  c \theta^\al \Big(F_2(R) + D_*(R) \Big) +  C(K,\theta) \ep_2 +  \\ +
c(\theta) ~\Big(  A^{\frac 23}(R) E^{\frac 43}(R) +  A^{\frac 56}_*(R)  F^{\frac 16}(R) E_*(R) \Big)
\endgathered
$$
Taking into account  (\ref{C}) and the obvious inequality $F_2(R)\le A_*(R)$ we arrive at
$$
\gathered
F_2(2\theta R) + D_*(2\theta R)  \ \le \  c \theta^\al \mathcal E(R) +  c(K,\theta) +
 \\ + c(\theta) \Big(\mathcal E^{\frac 23}(R) \ep_2^{\frac 43} +
\mathcal E^{\frac {11}{12}}(R) ( C(K)+ R^{\frac \al{12}} M^{\frac 1{12}}) K\Big)
\endgathered
$$
Applying the Young inequality we get
$$
\gathered
F_2(2\theta R) + D_*(2\theta R)  \ \le \  \Big(\frac 14 + c \theta^\al \Big) \mathcal E(R) +  C(\theta,K) (1+
R^\al M)
\endgathered
$$

Gathering the estimates we obtain
$$
\gathered
\mathcal E(\theta R)  \le    \l(\frac 14 + c \theta^\al \r) \mathcal E(R) +  C(\theta,K) (1+ R^\al M).
\endgathered
$$
Fixing $\theta\in (0,\frac 12)$ so that
$$
c \theta^\al \leq \frac 14
$$
Hence
$$
\gathered
\mathcal E(\theta R)  \le    \frac 12\mathcal E(R) +  C(K) (1+
R^\al M).
\endgathered
$$
Next with the help of technic used in the proof of lemma \ref{ckn_lemma1} we obtain (\ref{Estimate_varE}).

\qed


\begin{lemma} \label{Estimate_A_A_D} Assume all conditions of Theorem \ref{Boundedness_A_A_D} hold and fix $R_0\in (0,1)$ so that
\begin{equation}
R_0^\al M \ \le \ 1.
\label{rho_0}
\end{equation}
Then for any $0<r<R\le R_0 $  the following estimates hold:
\begin{equation}
\gathered
 A(r) + A_*(r)    \le  c\left(\frac{r}{R} \right)^{\ga}
\Big( A( R ) + A_*( R) \Big)   +   \ep^{\frac 14}D(R) +  G(K, \ep)
\endgathered
\label{Statement_of_Theorem_A_A}
\end{equation}
\begin{equation}
\gathered
D(r) \ \le \ c\left( \frac{r}{R}\right)^{\ga} D(R)\ + \ c(K) \Big(A^{\frac {11}{12}}(R)+A_*^{\frac {11}{12}}(R)\Big) + G(K, \ep)
\endgathered
\label{Cor}
\end{equation}
where $\ga>0$ is some absolute constant and $G$ is a continuous function possessing the following property:
\begin{equation}
\gathered
\mbox{for any fixed}\quad K >0 \quad G(K, \ep)\to 0 \quad \mbox{as}\quad\ep\to 0.
\endgathered
\label{Property of G}
\end{equation}
\end{lemma}

\noindent
{\bf Proof.} From (\ref{C_3}) taking into account (\ref{rho_0}) we obtain
\begin{equation}
C(r) \ \le A^{\frac 12}(r)\ep^{\frac 12}, \qquad F(r) \ \le \  A_*^{\frac 12}(r)C(K)
\label{F_K}
\end{equation}

Take arbitrary $r\in (0,R_0)$ and $\theta\in (0,\frac 12)$.
Denote by \ $\mathcal E_*(R)$ \ the following functional
$$
\mathcal E_* (R) \ = \ A(R)+   A_*(R),
$$
Then from  (\ref{Lo_En_Iq}) similar to (\ref{Similar to this}) using (\ref{F_K}) we derive
$$
\gathered
\mathcal E_*(\theta R)  \ \le \  F_2(2\theta R) \ + \  C^{\frac 12}(2\theta R)D^{\frac 12}(2\theta R) \ +
\\
+\ c(\theta)\Big( \mathcal E_*^{\frac 12}(R)\ep^{\frac 12} + \mathcal E_*^{\frac 34}(R)\ep^{\frac 34} +
\mathcal E_*^{\frac 34}(R) K^{\frac 12}\ep^{\frac 14} + \mathcal E_*^{\frac 34}(R) C(K) \ep^{\frac 12}\Big)
\endgathered
$$
Applying the Young inequality and using (\ref{D}) we obtain
\begin{equation}
\gathered
\mathcal E_*(\theta R)  \ \le \  \frac 18 ~\mathcal E_*(R) \ + \  c(\theta)G(K,\ep)  \ +
\\
+ \ F_2(2\theta R) \ + \  C^{\frac 12}(2\theta R)D^{\frac 12}_*(2\theta R)
\endgathered
\label{1}
\end{equation}

From (\ref{F_2}) we conclude
\begin{equation}
F(2\theta R) \ \le c\theta^\al  \mathcal E_*(R) \ + \ G(K, \ep).
\label{2}
\end{equation}

From (\ref{D_*}) for $R\le R_0$ with the help of (\ref{F_K}) and the Young inequality we obtain
\begin{equation}
\gathered
D_*(2\theta R) \ \le \ c ~\theta^{\be} D_*(R) \ + \ C(\theta,K)\mathcal E^{\frac {11}{12}}_*(R)  + c(\theta)G(K,\ep)
\endgathered
\label{DDD}
\end{equation}
Hence from (\ref{F_K}) we obtain
$$
\gathered
 C^{\frac 12}(2\theta R)D^{\frac 12}(2\theta R) \ \le \ c(\theta) \mathcal E_*^{\frac 14}(R)\ep^{\frac 14} D_*^{\frac 12}(R)    \ +  \\
 + \ C(\theta,K) \ep^{\frac 14} \mathcal E^{\frac {17}{24}}_*(R) + c(\theta)G(K,\ep)
\endgathered
$$
Applying the Young inequality we arrive at
\begin{equation}
\gathered
 C^{\frac 12}(2\theta R)D^{\frac 12}(2\theta R) \ \le \ \frac 18 \mathcal E_*(R)   +  \frac 12 \ep^{\frac 14} D_*(R)  + \ c(\theta)G(K, \ep)
\endgathered
\label{3}
\end{equation}

Gathering estimates (\ref{1}) --- (\ref{3}) we obtain the inequality
$$
\gathered
\mathcal E_*(\theta R)  \ \le \  \Big(\frac 14 + c\theta^\ga \Big) ~\mathcal E_*(R)  + \frac 12 \ep^{\frac 14} D_*(R) +  c(\theta)G(K,\ep)
\endgathered
$$
Choosing $\theta\in (0,\frac 12)$ so that
 $$
 \frac 14 + c\theta^\al \ = \  \frac 12
 $$
we obtain
$$
\gathered
\mathcal E_*(\theta R)  \ \le \  \frac 12 ~\mathcal E_*(R)  + \frac 12 \ep^{\frac 14} D_*(R) +  c(\theta)G(K,\ep)
\endgathered
$$
Iterating this inequality we obtain (\ref{Statement_of_Theorem_A_A}).

Choosing in (\ref{DDD}) $\theta\in (0,\frac 12)$ so that
$$
c\theta^\be \ =  \ \frac 12
$$
and iteration the inequality we derive (\ref{Cor}).

\qed


\begin{theorem}  \label{CKN_condition} For any $K>0$ there exists a constant $\ep_0(K)>0$ such that if the condition
(\ref{Statement_of_Theorem_E_E1}) holds with $\ep\le \ep_0$,  then
there exists $\rho_*\in (0,1)$ such that
$$
\begin{array}c
 \Big( C(\rho_*) + F(\rho_*) + D(\rho_*)\Big) \ < \ \ep_*^{\frac 13},
\end{array}
$$
where the constant $\ep_*>0$ is defined in Theorem \ref{Fixed_r}.
\end{theorem}

\noindent
{\bf Proof.} From (\ref{Estimate_varE}) we obtain
$$
\limsup\limits_{r\to 0} D_*(r) \le c(K).
$$
From (\ref{Statement_of_Theorem_A_A}) we derive
$$
\gathered
\limsup\limits_{r\to 0} \Big(A(r)+A_*(r)\Big) \ \le
\ \ep^{\frac 14} \limsup\limits_{\rho\to 0} D(\rho) + G(K, \ep ) \ \le \\ \le  \ \ep^{\frac 14} c(K) + G(K, \ep).
\endgathered
$$
From (\ref{Cor}) we obtain
$$
\gathered
\limsup\limits_{r\to 0} D_*(r) \ \le  \ c(K)\limsup\limits_{\rho\to 0} \Big(A^{\frac {11}{12}}(\rho)+ A^{\frac {11}{12}}_*(\rho)\Big)  + G(K, \ep ) \ \le \\
\le  \ c(K) \Big(\ep^{\frac 14} c(K) + G(K, \ep)\Big)^{\frac {11}{12}} + G(K, \ep).
\endgathered
$$
From (\ref{C_3}) we conclude
$$
\gathered
\limsup\limits_{r\to 0}  \Big( C(r)+F(r)\Big) \ \le (\ep^{\frac 12} + C(K)) \limsup\limits_{r\to 0} \Big(A(r)+A_*(r)\Big)  \ \le \\ \le \
(\ep^{\frac 12} + C(K))\Big( \ep^{\frac 14} c(K) + G(K, \ep)\Big)^{\frac 12}.
\endgathered
$$

Taking into account  (\ref{Property of G}) for any $K>0$ we can find $\ep_0(K)>0$ such that for any $\ep\in (0,\ep_0)$
$$
c(K) \Big(\ep^{\frac 14} c(K) + G(K, \ep)\Big)^{\frac {11}{12}}
 + G(K, \ep) \ < \ \frac{\ep_*^{\frac 13}}{2}
$$
and
$$
(\ep^{\frac 12} + K^{\frac 12})\Big( \ep^{\frac 14} c(K) + G(K, \ep)\Big)^{\frac 12} \ < \ \frac{\ep_*^{\frac 13}}{2}.
$$
Hence for $\ep\in (0,\ep_0)$
$$
\limsup\limits_{r\to 0}  \Big( C(r)+F(r)+ D_*(r)\Big) \ < \ \ep_*^{\frac 13}.
$$

\qed

\end{document}